\documentclass[english]{IEEEtran}

\usepackage{graphicx}
\usepackage{acronym}
\usepackage{amsmath}
 \usepackage{url,doi} 

\usepackage{cite}
\usepackage{amssymb}

\usepackage{balance}

\usepackage[ruled,linesnumbered,norelsize]{algorithm2e}

\makeatletter
\renewcommand{\@algocf@capt@plain}{above}
\makeatother
\usepackage[T1]{fontenc}

\newcommand{\measdim}{d} 
\newcommand{\statedim}{n} 
\newcommand{\jacobian}{\Psi} 
\newcommand{\meascov}{\Phi} 
\newcommand{\nonlvar}{\Upsilon} 
\newcommand{\mean}{\mu}

\newcommand{\priormean}{\mean^-}
\newcommand{\posteriormean}{\mean^+}

\newcommand{\meas}{y}
\newcommand{\transmeas}{\tilde y}
\newcommand{\predmeas}{\hat y}

\newcommand{\measfun}{h}
\newcommand{\transmeasfun}{\tilde \measfun}

\newcommand{\statecov}{P} 
\newcommand{\posteriorcov}{{\statecov^+}} 
\newcommand{\priorcov}{{\statecov^-}} 

\newcommand{\noisecov}{R} 

\newcommand{\matr}[1]{\begin{bmatrix}#1\end{bmatrix}}

\newcommand{\state}{x}

\newcommand{\measnoise}{\varepsilon}
\newcommand{\transmeasnoise}{\tilde \measnoise}
\newcommand{\transmat}{D}
\newcommand{\nonl}{\eta}
\newcommand{\nonllimit}{\eta_\text{limit}}
\newcommand{\innocov}{S}
\newcommand{\kalmangain}{K}

\newcommand{\eigVECTOR}{U}
\newcommand{\eigVALUE}{\Lambda}

\DeclareMathOperator{\kld}{KLD}

\newcommand{\KLD}[1]{\kld\left(#1\right)}
\newcommand{\msqrt}[1]{\sqrt{#1}}

\newcommand{\el}[1]{_{\left[#1\right]}}

\newtheorem{example}{Example}

\newtheorem{proof}{Proof}

\DeclareMathOperator{\cov}{cov}
\DeclareMathOperator{\tr}{tr}
\DeclareMathOperator{\N}{N}

\renewcommand{\d}{\mathrm{d}}

\SetCommentSty{mycommfont}

\begin{document}
\newacro{BinoGMF}{binomial Gaussian mixture filter}
\newacro{CDF}{cumulative distribution function}
\newacro{KF}{Kalman filter}
\newacro{CKF}{cubature Kalman filter}
\newacro{EKF}{extended Kalman filter}
\newacro{nEKF2}{numerical second order Extended Kalman filter}
\newacro{EKF2}{second order Extended Kalman filter}
\newacro{IPLF}{iterated posterior linearisation filter}
\newacro{GGF}{general Gaussian filter}
\newacro{GMF}{Gaussian mixture filter}
\newacro{IKF}[IEKF]{iterated extended Kalman filter}
\newacro{IPLF}{iterated posterior linearization filter}

\newacro{KL}{Kullback-Leibler}
\newacro{KLD}{Kullback-Leibler divergence}
\newacro{KLPUKF}{Kullback-Leibler partitioned update Kalman filter}
\newacro{PUKF}{partitioned update Kalman filter}
\newacro{PDF}{probability density function}
\newacro{PF}{particle filter}
\newacro{SLF}{statistically linearized filter}

\newacro{S2KF}[S\textsuperscript{2}KF]{smart sampling Kalman filter}

\newacro{RUF}{recursive update filter}
\newacro{UKF}{unscented Kalman filter}

\title{Kullback-Leibler Divergence Approach to Partitioned Update Kalman Filter}
\author{Matti Raitoharju,  \'Angel F. Garc\'ia-Fern\'andez, and Robert Pich\'e}
\maketitle
\allowdisplaybreaks

\begin{abstract}
Kalman filtering is a widely used framework for Bayesian estimation. The partitioned update Kalman filter applies a Kalman filter update in parts so that the most linear parts of measurements are applied first. In this paper, we generalize partitioned update Kalman filter, which requires the use oft the second order extended Kalman filter, so that it can be used with any Kalman filter extension. To do so, we use a Kullback-Leibler divergence approach to measure the nonlinearity of the measurement, which is theoretically more sound than the nonlinearity measure used in the original partitioned update Kalman filter. Results show that the use of the proposed partitioned update filter improves the estimation accuracy.
\end{abstract}

\begin{IEEEkeywords}Bayesian estimation; nonlinear; estimation; Kalman filters; Kullback-Leibler divergence\end{IEEEkeywords}

\section{Introduction}
In Bayesian filtering we are interested in calculating the \ac{PDF} of a dynamic state based on a sequence of measurements. It is a recursive process in which a prior distribution is updated using a measurement to obtain a posterior distribution. This distribution then evolves in time to become a new prior distribution. Bayesian estimation has a wide range of applications from positioning \cite{6746180}, tracking \cite{Sadhu20063769}, and quality control \cite{Zhai20102319}, to brain imaging \cite{hiltunen2011state} and modeling spread of infectious diseases \cite{opac-b1124487}. In general, a Bayesian estimate cannot be computed in closed form. Under certain conditions, which include that measurements and state transition function have to be linear and associated noises Gaussian, Bayesian estimates can be computed in closed form using the algorithm known as the  \cite{kalman}. For nonlinear measurements various Kalman filter extensions have been developed. When the nonlinearity is small, the \ac{KF} extensions produce accurate estimates, but large nonlinearity can cause serious inaccuracies. 

Several methods have been proposed to quantify the amount of nonlinearity in order to monitor the performance of a \ac{KF} extension \cite{jazwinski,raitoharju2014a,splitmerge,HavlakandCampbell,Huber,6916118,6584787}. An algorithm based on~\ac{EKF2} was presented in \cite[p. 349]{jazwinski} and it was extended for multidimensional correlated measurements in \cite{raitoharju2014a}. Although it was found to be a good indicator of the accuracy of a \ac{KF} extension compared to methods that are computationally feasible presented in~\cite{splitmerge,HavlakandCampbell,Huber}, its evaluation requires computation of second order derivatives of the measurement function or computation of the measurement function values in a number of points that increases quadratically.  

The effect of the nonlinearity can also be gauged by comparing the higher order moments of the posterior to those corresponding to the Gaussian assumption. This is considered for example in~\cite{6916118}.

If the prior covariance is sufficiently small, any nonlinear function is well approximated as a linear function in the prior's high-probability region. Therefore, if a nonlinearity measure indicates that a measurement is highly nonlinear, the prior can be split into a sum of small-covariance parts that can be used as the components in a \ac{GMF} \cite{Sorenson_Alspach_1971}.  As the computational burden increases with every split, the number of new components needed to bring the nonlinearity measure below a certain threshold should be minimized \cite{BINOGMF}.

In~\cite{PUKF} another use of the  \ac{EKF2}-based  nonlinearity measure  in~\cite{raitoharju2014a} was presented. The algorithm, called the \ac{PUKF}, transforms a measurement vector with a linear transformation so that the nonlinearity of the least nonlinear measurement element is minimized. Only the elements with nonlinearity below a threshold are applied first and then the linearization and evaluation of the nonlinearity is redone. This way the estimates become more accurate.

In~\cite{6584787}, it was shown that the moments computed in a \ac{GGF} \cite{855552} can be used to compute the  \ac{KLD} \cite{Kullback51klDivergence} of the \ac{GGF} approximation to the joint \ac{PDF} of the state and the measurement from the true joint \ac{PDF}. Various \ac{KF} extensions can be seen as approximations of the \ac{GGF} and, thus, they can be used to approximate this \ac{KLD}, which can be also interpreted as a nonlinearity measure. This \ac{KLD} has been used in Kalman optimization in \cite{6963359}.

In this paper, we present an algorithm that is similar to \ac{PUKF}, but uses the \ac{KLD} based nonlinearity measure. The benefits of using the \ac{KLD} measure is that it is mathematically sound and it can be used in combination with any \ac{GGF} approximation that can be used to approximate the \ac{KLD}.  If the proposed algorithm is used with \ac{EKF2} or a numerical approximation of it, the results are the same as with the original \ac{PUKF}.

The rest of this paper is organized as follows: Section~\ref{sec:background} gives the background work. The new algorithm is developed in Section~\ref{sec:KLPUKF}. Section~\ref{sec:examples} presents examples of the use and accuracy of the proposed algorithm. Section~\ref{sec:conclusions} concludes the paper.

\section{Background work}
\label{sec:background}
In this paper, we  derive an improved version of the \ac{PUKF}. To do so, it is convenient to introduce some background material. In Section~\ref{sec:formulation}, we first present the \ac{GGF} update and a measure of its performance based on the \ac{KLD}. In Section~\ref{sec:PUKF}, we revisit the PUKF. 

\subsection{\ac{GGF}}
\label{sec:formulation}
In this paper, we consider the Bayesian update step of a state. We assume that state $x \in \mathbb{R}^{\statedim}$ has a Gaussian prior \ac{PDF}
\begin{equation}
p(\state)=\mathrm{N}(\state | \priormean, \priorcov).
\end{equation}
where $\mathrm{N}$  is the \ac{PDF} of a normal distribution, $\priormean$ is the prior mean, and $\priorcov$ is the prior covariance. This state is observed through a measurement that is modeled with a measurement model of form
\begin{equation}
	\meas= \measfun( \state) + \measnoise \label{equ:meg},
\end{equation}
where $\meas$ is the $\measdim$-dimensional measurement value and $\measnoise$ is a zero-mean Gaussian measurement noise with covariance $\noisecov$.

The objective is to compute the posterior 
\begin{equation}
	p( x | y) \propto p(x)p(y|x), \label{equ:bayes}
\end{equation}
where $\propto$ stands for proportionality and $p(y|x) = \N (y| \measfun(x), R ) $ is the density of the measurement given the state. If function $\measfun(\cdot)$ is linear, the posterior can be computed exactly using the \ac{KF}. But when   $\measfun(\cdot)$ is not linear, an approximation has to be used.

A general way to formulate a Kalman filter type approximation is to use the \ac{GGF} formulation for additive noise \cite{855552}. 
In the update step, we approximate the posterior as a Gaussian with mean $\posteriormean$ and covariance $\posteriorcov$ by
\begin{align}
	\posteriormean & = \priormean + \kalmangain(\meas - \predmeas) \label{equ:mean}\\
	\posteriorcov & = \priorcov - \kalmangain\innocov\kalmangain^T, \label{equ:posteriorcov}
\end{align}
where
\begin{align}
	\predmeas &= \int \measfun(\state)p(\state) \d \state \label{equ:i1}\\
	\jacobian &= \int \left( \state - \priormean\right) \left( \measfun(\state)  - \predmeas\right)^Tp(\state) \d \state \label{equ:i2} \\
	\meascov &= \int \left( \measfun(\state)  - \predmeas\right)\left( \measfun(\state)  - \predmeas\right)^Tp(\state)  \d \state \label{equ:i3} \\
	\innocov & =\meascov + \noisecov \label{equ:innocov} \\ 
	\kalmangain &= \jacobian \innocov^{-1} \label{equ:kgain}.
\end{align}

The \ac{GGF} is implicitly defining a Gaussian approximation of the joint density of state $\state$ and measurement $\meas$. The approximation, denoted $q(x,y)$ is not the exact joint density $p(x,y)$. The approximation error can be measured using \ac{KLD} \cite{Kullback51klDivergence}
\begin{align}
\eta= \KLD{ p,q } &= \int{\int p(x,y) \log{\frac{p(x,y)}{q(x,y)}} \d x}\d y
\end{align}
In \cite{6584787}, it was shown that 
\begin{equation}
	\nonl =\frac{1}{2}\log \left| I + \noisecov^{-1}\nonlvar  \right|, \label{equ:nonl}
\end{equation}
where
\begin{equation}
	\nonlvar  = \meascov - \jacobian^T \left(\priorcov\right)^{-1} \jacobian . \label{equ:nonn}
\end{equation}
Different \ac{KF} extensions,  such as  the \ac{UKF} \cite{WANUKF} and the \ac{CKF} \cite{cubature}, can be interpreted as  approximations of the \ac{GGF} and they can also be used to compute approximations of the integrals \eqref{equ:i1}-\eqref{equ:i3}  and, thus, approximate the joint \ac{KLD}. Some \ac{KF} extensions, such as the \ac{EKF}, make linearizations such that the moments \eqref{equ:i2} and \eqref{equ:i3} can be written in the form
\begin{align}
	\jacobian & = \priorcov A^T \\
	\meascov & = A \priorcov A^T,
\end{align}
where $A$ is a $\measdim \times \statedim$ matrix. For these filters $\nonlvar$ in \eqref{equ:nonn} is always 0 and they cannot be used to approximate the \ac{KLD}.

\subsection{\ac{PUKF}}
\label{sec:PUKF}

For conditionally independent measurements $\meas_1,\ldots, \meas_\measdim$ given the state we can write
\begin{align}
	p(\meas_{1:\measdim} | \state)=\prod_{i=1}^\measdim p(\meas_i | \state). \label{equ:condind}
\end{align}
The prior can be updated using the measurements sequentially:
\begin{equation}
\begin{aligned}
	p( x | \meas_{1:n}) &\propto  p(\state)p(\meas_1 | \state)  \prod_{i=2}^\measdim p(\meas_i | \state) \\
	 &\propto p(\state | \meas_{1})  \prod_{i=2}^\measdim p(\meas_i | \state).
\end{aligned}
\end{equation}
When the measurement model is of the form \eqref{equ:meg} and the noise covariance $\noisecov$ is diagonal, the measurements are conditionally independent.
Thus, in a linear Gaussian measurement model, the \ac{KF} update can be applied one measurement element at a time \cite[p. 119]{Maybeck79} and the posterior distribution does not change.  However, when the updates are approximate, the final posterior \ac{PDF} approximation changes.

In general,  the $\measdim$ measurement elements are not conditionally independent given the state so \eqref{equ:condind} does not hold. The main idea behind the \ac{PUKF} \cite{PUKF} is to apply a linear transformation to the measurement model so that measurement elements are conditionally independent and the measurements are applied in an order that aims to minimize the approximation error.  A general description of \ac{PUKF} is given in Algorithm~\ref{algo:pukf}.
 \begin{algorithm}[b]
 \caption{A general description of the \ac{PUKF} algorithm}
 \label{algo:pukf}
 \SetNlSty{textbf}{\small}{}
Evaluate the nonlinearity of the measurement elements.\\
Minimize the nonlinearity of part of the measurement by applying a linear transformation to the measurement. \\
 Update the state using the part of the measurement whose nonlinearity is smaller than a set threshold. \\
 If the whole measurement is not applied, use the partially updated prior as a new prior and the unused measurements as a new measurement and return to 1. \\
 Return posterior. \\
\end{algorithm}

The \ac{PUKF} in \cite{PUKF} can only be used with \ac{EKF2} \cite[pp. 345-347]{jazwinski}, which is based on second order Taylor expansion of the measurement function, or with a central difference filter  \cite{855552} that is a numerical approximation of \ac{EKF2}. We call this numerical approximation the \ac{nEKF2}. The \ac{EKF2}-specific nonlinearity measure is
\begin{equation}
\hat \eta = \sum_{i=1}^{\measdim} \sum_{j=1}^{\measdim}  \left( R^{-1} \right)\el{i,j}\tr PH_iP H_j,
\end{equation}
where $H_i$ is the Hessian of $i$th element of the measurement function $h(\cdot)$ or its numerical approximation and subscript $[i,j]$ means the matrix element in the $i$th row and $j$th column.

\begin{example}
We proceed to illustrate how the \ac{PUKF} works as it is the foundation of the methods we propose in this paper. We consider an example from \cite{PUKF}. Here the prior is one dimensional with \ac{PDF} $\N(x|1,1)$ and the nonlinear measurement equation is \begin{equation}
y=\matr{x^2 - 2x - 4 \\
   -x^2 + \frac{3}{2}} + \varepsilon,
\end{equation}
where $\varepsilon$ has \ac{PDF} $\N(\varepsilon|0,I)$. The measurement can be transformed with the linear transformation
\begin{equation}
	\hat{y} = \frac{1}{\sqrt{2}}\matr{1 & 1 \\ 1 &-1} y
\end{equation}
into a linear term and a polynomial term:
\begin{equation}
\hat y= \sqrt{2} \matr{-x - \frac{5}{4}   \\
   x^2 - x - \frac{11}{4}} + \hat\varepsilon,
\end{equation}
where $\hat\varepsilon \sim \N(0,I)$. In \ac{PUKF} the linear measurement function is applied first and the partially updated state has mean $-\frac{1}{2}$ and covariance $\frac{1}{3}$. The polynomial measurement function is applied using this partially updated state and \ac{EKF2}-based update. In Figure~\ref{fig:pdfs} \ac{EKF2} is used as a reference. \ac{EKF2} applies both measurements at once and the posterior estimate is the same for the original and transformed measurement models. When compared to the true posterior, which is computed using a dense grid, the posterior estimate of \ac{PUKF} is significantly more accurate than the \ac{EKF2} posterior estimate.
\begin{figure}
	\centering
	\includegraphics[width=0.8\columnwidth,clip=true,trim=1cm 0.5cm 1.2cm 0.5cm]{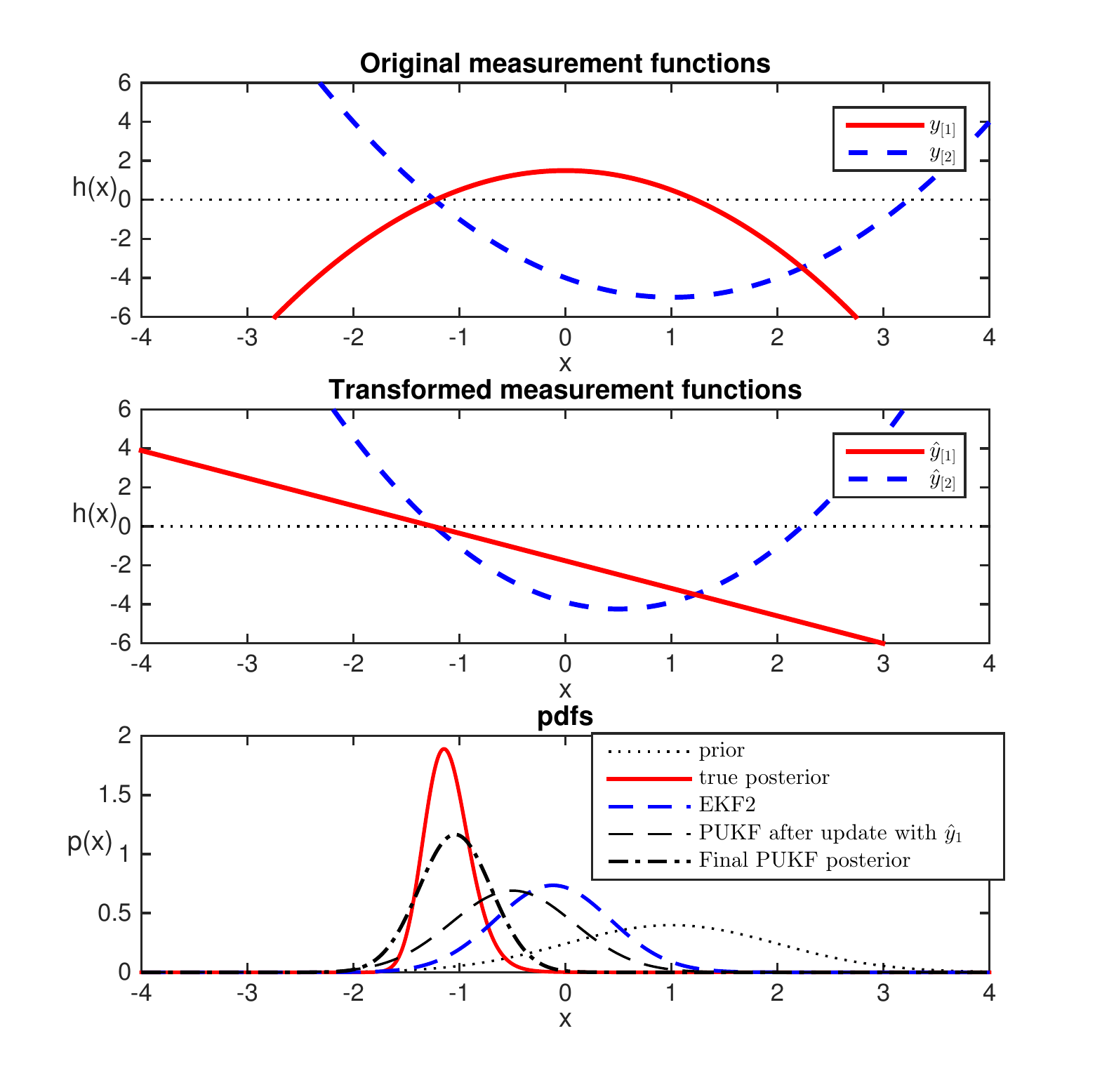}
	\caption{Transforming second order polynomial measurements to minimize nonlinearity of $\hat{\meas}_1$ and posterior comparison of  \ac{PUKF} and \ac{EKF2}}
	\label{fig:pdfs}
\end{figure}
\end{example}

\section{\ac{KLPUKF}}
\label{sec:KLPUKF}
In this paper we develop \ac{KLPUKF} that uses \eqref{equ:nonl} as the nonlinearity measure. Thus, algorithm is not limited to \ac{EKF2} based filters and any filter, such as \ac{UKF} \cite{WANUKF} or \ac{CKF} \cite{cubature}, that is a good enough approximation of the \ac{GGF} can be used. As noted earlier, some algorithms, such as \ac{EKF}, produce always 0 nonlinearity. Thus, they are not good enough approximations of the \ac{GGF}. We assume that the \ac{KF} extension does the moment approximations at the prior. Iterative Kalman filter extensions, such as \ac{IKF} \cite[pp. 349-351]{jazwinski} and \ac{IPLF} \cite{PLF}, are not considered in this paper to be used with the \ac{KLPUKF}.

We recall that the steps of the \ac{PUKF} algorithm are given in  Algorithm~\ref{algo:pukf}. Therefore, in this section we specify how to compute the nonlinearity and the linear transformation. The linear transformation of the measurement model  can be written as
\begin{equation}
\transmeas= \transmat \meas = \transmat \measfun(\state) + \transmat \measnoise = \transmeasfun(\state) + \transmeasnoise \label{equ:transmeas},
\end{equation}
where $\transmat$ is a nonsingular square matrix.  In the appendix  we show that the linear transformation \eqref{equ:transmeas} does not change the posterior if all measurements are applied at once (Proof~\ref{proof:estc}) and does not have an effect on the nonlinearity measure \eqref{equ:nonl} (Proof~\ref{proof:nonldoesnotchange}). 

However, if an update is split into a sequence of independent measurement updates, the posterior may change. We use  a transformation such that the transformed noise covariance $\tilde R$ is
\begin{equation}
\tilde R=\cov \transmeasnoise=DRD^T = I. \label{equ:Icov}
\end{equation}
We show in Proof~\ref{proof:nonldoesnotchange} in the Appendix that the nonlinearity \eqref{equ:nonl} can be written now as
\begin{equation}
	\nonl =\frac{1}{2}\log \left| I + \tilde\nonlvar  \right|, 
\end{equation}
where
\begin{equation}
\tilde \nonlvar = D\nonlvar D^T. \label{equ:diagcov}
\end{equation}
We further want  $\tilde \nonlvar$ to be diagonal. Now the nonlinearity measure  of the $i$th measurement element can be defined as
\begin{equation}
\tilde \eta_i = \frac{1}{2}\log\left(1+ \tilde\nonlvar_{[i,i]}\right).
\end{equation}
When $\tilde R=I$ and $\tilde \nonlvar$ is diagonal the nonlinearity of the least nonlinear measurement element is minimized. This is shown in Proof~\ref{proof:maxnl} in the Appendix. Thus, this gives the best choice for a single element update.

Rather than choosing $\tilde{R}=I$, we could have chosen any other diagonal matrix. However, this would increase the value of the corresponding diagonal values of $\tilde{\nonlvar}$ and the \ac{KLD} values of individual measurements of \eqref{equ:nonlvalue} would not change. Thus, choosing $\tilde{R}=I$ is an arbitrary, but natural choice.

To compute the matrix $D$ we first introduce the notation for matrix square root
\begin{equation}
	\msqrt{\noisecov}\msqrt{\noisecov}^T=\noisecov, \label{equ:sqrtr}
\end{equation}
which can be computed with, e.g., Cholesky decomposition. The transformation $\transmat$ that makes $\tilde{\nonlvar}$ diagonal and $\tilde \noisecov=I$ is \cite{PUKF}
\begin{equation}
	\transmat = \eigVECTOR^{T} \msqrt{\noisecov}^{-1}, \label{equ:D}
\end{equation}
where $\eigVECTOR$ is computed using an eigendecomposition
\begin{equation}
	\eigVECTOR\eigVALUE\eigVECTOR^T = \msqrt{\noisecov}^{-1} \nonlvar \msqrt{\noisecov}^{-T}, \label{equ:eigdecom}
\end{equation}
where $\eigVECTOR$ is orthogonal and the eigenvalues in the diagonal matrix $\eigVALUE$ are sorted in ascending order.
We can see that this transformation fulfills our requirements of having identity transformed measurement noise covariance \eqref{equ:Icov} from 
\begin{equation}
	\tilde{R} = \eigVECTOR^T \msqrt{\noisecov}^{-1} R \msqrt{\noisecov}^{-T} \eigVECTOR = I
\end{equation}
and the diagonality of matrix $\tilde \noisecov$ \eqref{equ:diagcov} can be seen from orthogonality and \eqref{equ:eigdecom}
\begin{equation}
	\tilde{\Upsilon} = \eigVECTOR^T \msqrt{\noisecov}^{-1} \Upsilon \msqrt{\noisecov}^{-T} \eigVECTOR = \eigVALUE.
\end{equation}

To save computational resources, rather than using only one measurement in each update, we perform the update using all measurement elements that have $\tilde \nonl_i$ below a limit $\eta_\text{limit}$ or if this set is empty then with the measurement element with smallest $\tilde \nonl_i$. We also note that after a measurement update the moments \eqref{equ:i1}-\eqref{equ:i3} are recomputed. Because the nonlinearities of the measurements  change, the transformation matrix $\transmat$ is recomputed between every update.

The \ac{KLPUKF} algorithm is given in Algorithm~\ref{algo:KLPUKF}. When using \ac{EKF2} or a numerical approximation of it as basis of the filter and setting $\nonllimit$ properly, the algorithm produces the same results as with the original \ac{PUKF}.

\begin{algorithm}
\caption{Measurement update step in \ac{KLPUKF}}
\label{algo:KLPUKF}
 \SetNlSty{textbf}{\small}{}
       \SetKwInOut{Input}{input}
        \SetKwInOut{Output}{output}
        \SetKwInOut{Parameters}{parameters}
        \Input{
        Prior state: $\mean_0$ --  mean 
                  $\statecov_0$ --  covariance \newline
        Measurement model: $\meas_0$ --  value, $\measfun_0(\cdot)$ --  function, 
                  $\noisecov_0$~--~covariance \newline
                  }
         \Output{
         Updated state: $\mean^+$ -- mean, $\statecov^+$ -- covariance
         }  
 Compute $\msqrt{R}$ (\ref{equ:sqrtr})\\
 $d \leftarrow \text{initial measurement dimension}$ \\
 $i \leftarrow 0$ \tcp{Iteration counter}
 \While{ $\measdim > 0$}
 {
 Compute approximations  ($\predmeas_i$, $\jacobian_i$, and $\meascov_i$)  of moments \eqref{equ:i1}-\eqref{equ:i3} for the $\measdim$-dimensional measurement $h_i(\cdot)$ using a \ac{KF} extension and prior $\N(x| \mean_i,\statecov_i)$\\
 $\nonlvar  \leftarrow \meascov - \jacobian^T \left(\statecov_i\right)^{-1} \jacobian$\\
   Compute $\eigVECTOR_i$ and $\eigVALUE_i$  from $\msqrt{R}$  and $\nonlvar$ using (\ref{equ:eigdecom}) \\
 $\transmat_i \leftarrow \eigVECTOR_i^T\msqrt{\noisecov_i}^{-1} $ \\
 Choose largest $k$ so that $\log( 1+ \eigVALUE_{i,[k,k]})\leq\nonllimit$. If no such $k$ exists, set $k\leftarrow1$\\
 
 	\tcp{Compute partial update}
	\Indp
 	$\tilde \predmeas_i \leftarrow {\transmat_i}\el{1:k,:} \predmeas_i   $ \\
	$\tilde\innocov_i \leftarrow {\transmat_i}\el{1:k,:}\meascov\transmat\el{1:k,:}^T+I$	\\
	$\tilde\kalmangain_i \leftarrow \jacobian_i^T{\transmat_i}\el{1:k,:}^T \tilde{\innocov}_i^{-1}$ \\
	$\mean_{i+1} \leftarrow \mean_{i} + \tilde{\kalmangain}_{i} ( {\transmat_i}\el{1:k,:}\meas_i - \tilde{\predmeas})$ \\
	$\statecov_{i+1} \leftarrow \statecov_{i} - \tilde{\kalmangain}_{i} \tilde{\innocov}_i \tilde{\kalmangain}_i^{T}$ \\
	\Indm
 	\tcp{Update remaining measurement}
	\Indp
	$\meas_{i+1} \leftarrow  {\transmat_i}\el{k+1:\measdim,:}\meas$ \\
	$\measfun_{i+1}( \state) \leftarrow {\transmat_i}\el{k+1:\measdim,:} \measfun_i(\state)$ \\
	$\msqrt{R}_{i+1} \leftarrow I$ \tcp{Updated measurement noise covariance is an identity matrix due to decorrelation}
	$d\leftarrow d-k$ \tcp{Updated measurement dimension}
	$i \leftarrow i+1$  \tcp{Increase counter}
	\Indm
}
$\mean^+\leftarrow \mean_i$  \tcp{Posterior mean}
$\statecov^+\leftarrow \statecov_i$  \tcp{Posterior covariance}
\end{algorithm}

\section{Simulation examples}
\label{sec:examples}
In this section, we present examples of how the \ac{KLPUKF} improves estimation accuracy. First, we evaluate the accuracy enhancements for a one-dimensional state that is observed through measurements containing trigonometric functions. Then we present examples with range measurements to beacons. In these examples we also test filtering and evaluate the effect of the nonlinearity limit $\nonllimit$. In our examples, we use an approximation of \ac{GGF} computed in a dense grid, \ac{UKF} \cite{WANUKF}, \ac{nEKF2} \cite{855552}, and \ac{S2KF} \cite{steinbring2013s}.

\begin{example} 1D example

First we consider an example with a 1-dimensional state and standard normal prior. The prior is updated with a 3-element measurement of form
\begin{equation}
y=\matr{ x+4 \sin(x)+7 \\
-x+4\sin(x)-4 \\
 -2\cos(x)-8} + \measnoise,
\end{equation}
where \ac{PDF} of $\measnoise$ is $\N(\measnoise| 0,I)$. The top-left plot of Figure~\ref{fig:sincosexample} shows these measurement functions. The dashed lines show the statistical linearizations of the functions obtained from moments \eqref{equ:i1}-\eqref{equ:i3} \cite{PLF}. The dashed black line is the prior mean and dotted lines show the $2\sigma$ limits; i.e.\ 95\% of the probability is within these limits.
\begin{figure}
\centering
	\includegraphics[width=0.8\columnwidth,clip=true,trim=1.7cm 1.2cm 1.5cm 0.8cm]{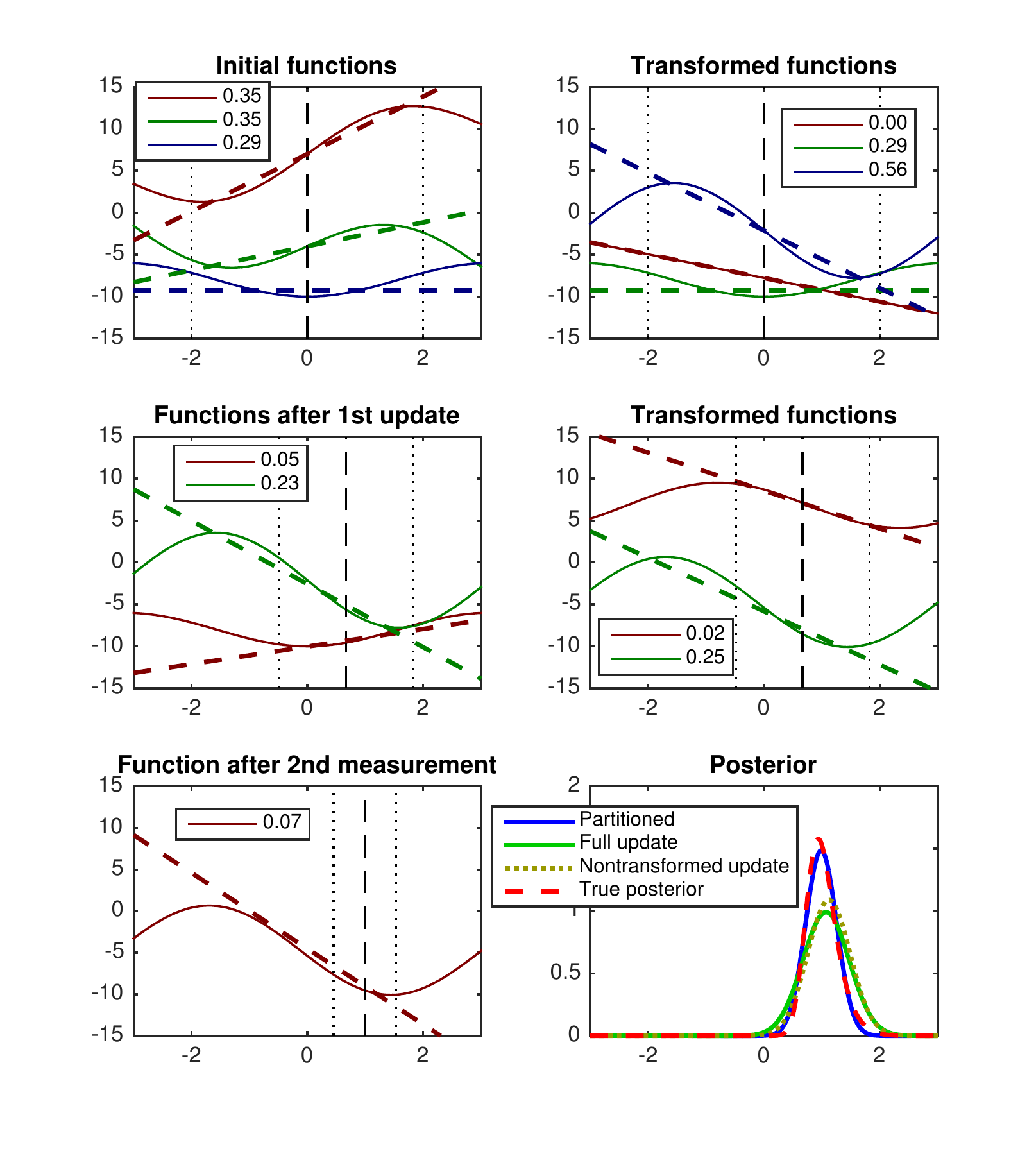}
	\caption{Example updates using all measurements at once and  using \ac{KLPUKF}. Red corresponds to $y_1$, green to $y_2$, and blue to $y_3$.}
	\label{fig:sincosexample}
\end{figure}
The numbers in the legend are the \acp{KLD} \eqref{equ:nonlvalue} of each element of the measurement function computed using \eqref{equ:nonlvalue}. Because $\nonlvar$ is not diagonal for initial functions, their \acp{KLD} does not sum to the total \ac{KLD} that is $0.8533$ for the initial measurements. The sum of \acp{KLD} of transformed functions is the same as the total \ac{KLD}.

Before transformation, the function with the cosine term (blue) has the smallest \ac{KLD}. For these measurements and prior the transformation matrix is
\begin{equation}
\transmat =\matr{ -\frac{1}{\sqrt{2}} & \frac{1}{\sqrt{2}} & 0 \\
		          0 & 0 & 1 \\
		          -\frac{1}{\sqrt{2}}  & -\frac{1}{\sqrt{2}} & 0},
  \end{equation} 
which causes the sin terms in the new first component to cancel out and the new measurements are
\begin{equation}
y=\matr{ -\sqrt{2}x - \frac{11}{\sqrt{2}}  \\
-2\cos x - 4 \\
  -4\sqrt{2} \sin(x) \frac{3}{\sqrt{2}}}  + \measnoise.
\end{equation}
The top-right plot shows these transformed measurements and corresponding \acp{KLD}. The linear measurement naturally has a zero \ac{KLD} and is applied first. The second row shows the remaining components of the transformed functions and the new linearizations that are made in the new prior. The last row shows the last measurement and the posterior \acp{PDF}. The curve labeled ``partitioned'' is the posterior computed with the proposed algorithm. ``Full update'' is the update using all measurements simultaneously and ``nontransformed update'' is the update made in parts without applying transformation $\transmat$ i.e.\ the measurement with the cosine element is applied first because it has the smallest \ac{KLD}. This plot shows that the posterior computed with the proposed method is the most accurate one.
\end{example}

\begin{example} Range measurements

In this example, we consider a more realistic situation with 2-dimensional state with prior \ac{PDF} $\N(x | \mathbf{0}, 12I)$ updated with 3 range measurements. The measurement model is 
\begin{equation}
	h(x)=\matr{\sqrt{ (x_1-2)^2 + (x_2-2)^2 } \\
	\sqrt{ (x_1+6)^2 + (x_2-6)^2 }\\
	\sqrt{ (x_1+2)^2 + (x_2-1)^2 }} + \measnoise, 
\end{equation}
where $\measnoise$ has \ac{PDF} $\N(\measnoise| 0,I)$ and the measurement values are $\matr{5 & 11.5 & 3.5}^T$.

Figure~\ref{fig:example1} shows an example of the update with \ac{KLPUKF}.  
The top row in Figure~\ref{fig:example1} shows the likelihoods of the range measurements (blue), prior (red) and the measurement likelihoods when linearized using the \ac{GGF} (orange). The second row shows the transformed measurements and their linearizations within the original prior. The last row shows the linearizations within the prior and how the prior and linearizations change after partial updates.  The integrals for \ac{GGF} were computed using a dense grid.
\begin{figure}
	\includegraphics[width=\columnwidth,clip=true,trim=2cm 0.8cm 1.6cm 0cm]{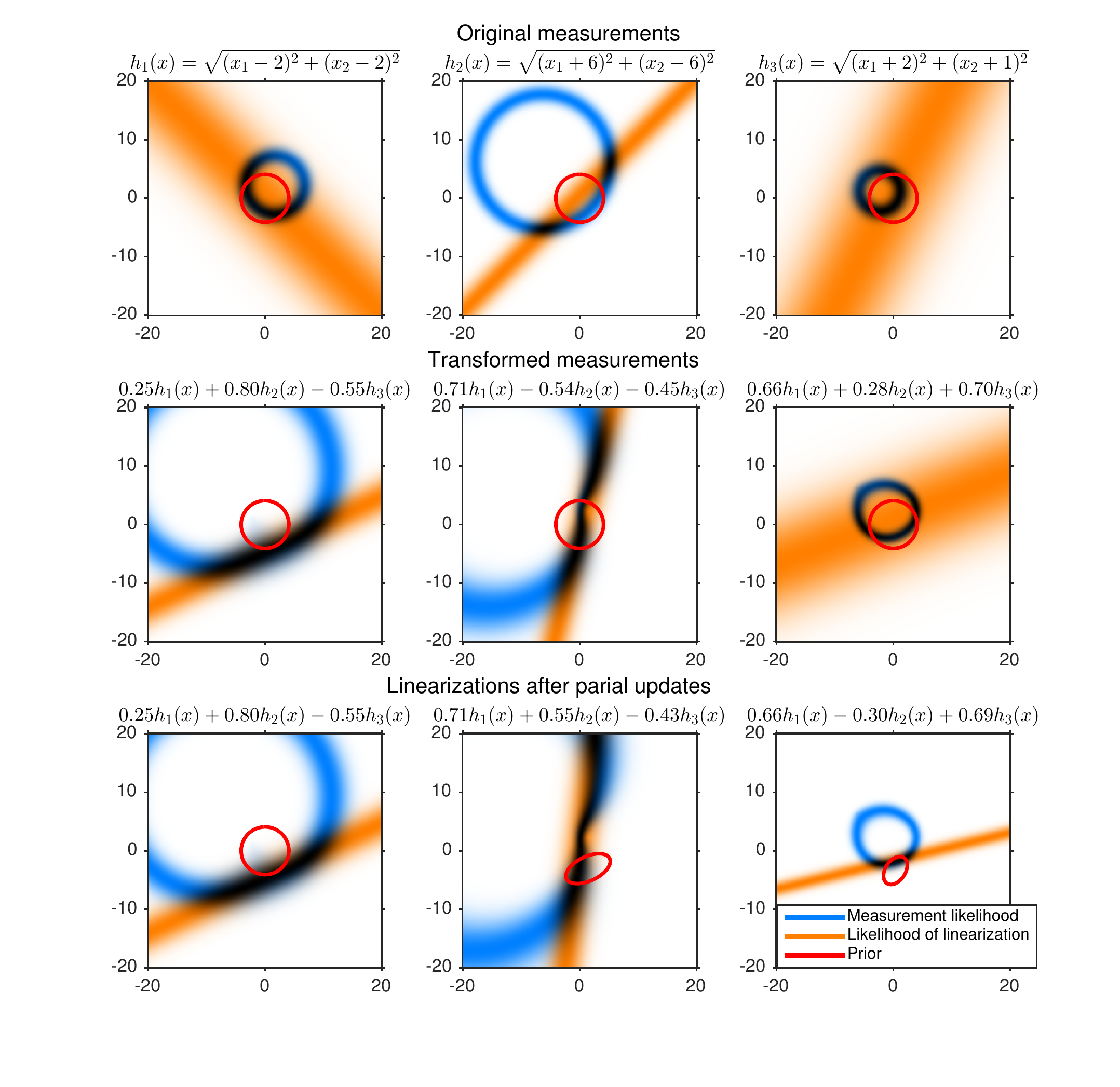}
	\caption{Example updates with using all measurements at once and with using \ac{KLPUKF} with \ac{GGF}. Ellipses contain 50\% of the prior probability. Intensity of the color of the likelihoods represents its value.}
	\label{fig:example1}
\end{figure}

We can see that the third measurement of both non-transformed and transformed measurements has a large uncertainty when linearized using the original prior. When using the partially updated prior the third measurement has much smaller uncertainty.

Figure~\ref{fig:example2} shows the same example computed using the \ac{UKF} linearizations.
\begin{figure}
	\includegraphics[width=\columnwidth,clip=true,trim=2cm 0.8cm 1.6cm 0cm]{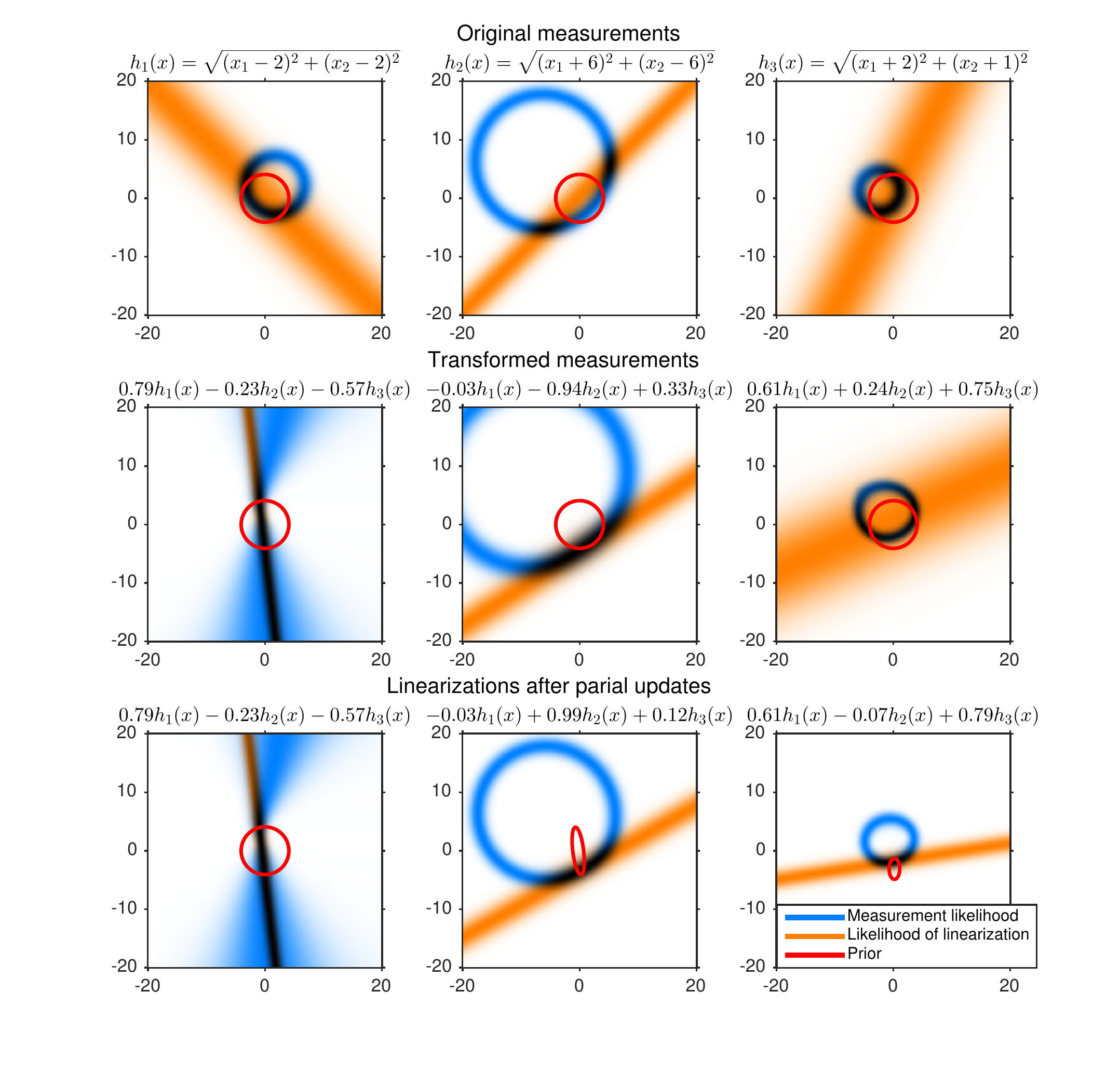}
	\caption{Example updates with using all measurements at once and with using \ac{KLPUKF} and \ac{UKF}.  Ellipses contain 50\% of the prior probability. Intensity of the color of the likelihoods represent its value.}
	\label{fig:example2}
\end{figure}
Contours containing 50\% of the posterior probability computed using \ac{GGF}, \ac{UKF}, \ac{nEKF2}  and their \ac{KLPUKF} counterparts are shown in Figure~\ref{fig:postres}. The contour for true posterior is computed in a grid. Figure shows how the \ac{GGF} produces a larger covariance than the true covariance and how the use of \ac{KLPUKF} makes the estimate closer to the true estimate. The \ac{UKF} estimate has a bad shape: it is too narrow horizontally and too long vertically. Using~\ac{KLPUKF} with \ac{UKF} reduces the covariance in the too long direction and the resulting posterior estimate is closer to the true posterior.
\begin{figure}
\centering
	\includegraphics[width=0.8\columnwidth]{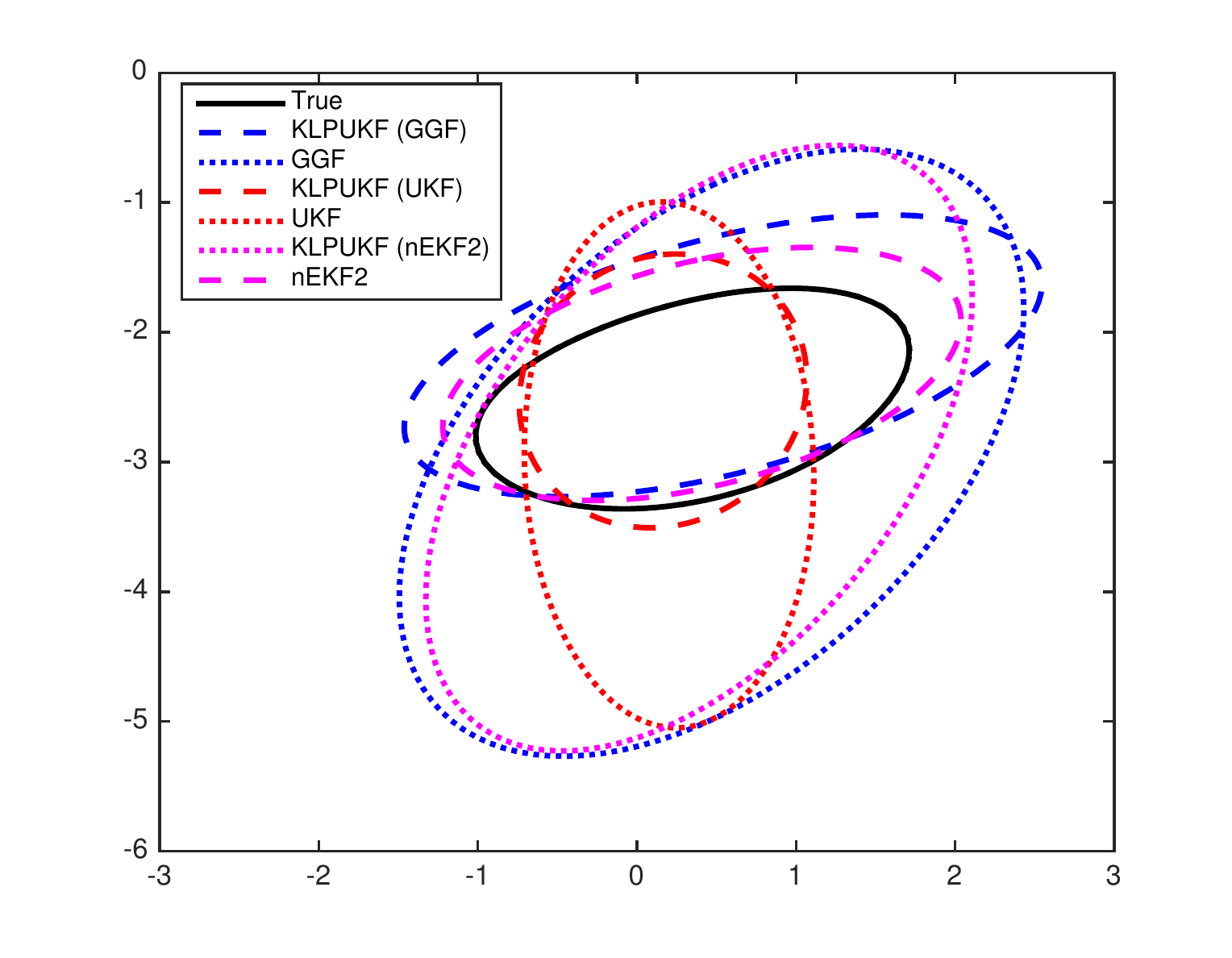}
	\caption{Contours containing 50\% of the posterior probability obtained with different algorithms}
	\label{fig:postres}
\end{figure}

\end{example}

\begin{example} Filtering example \nopagebreak

In the filtering example we consider the same measurement model as in the previous example. We use a 4-dimensional state 
\begin{equation}
	x=\matr{r_1 & r_2 & v_1 & v_2}^T,
\end{equation}
where $r_1$ and $r_2$ are position variables and $v_1$ and $v_2$ velocity variables. The state transition model is
\begin{equation}
	x_{t+1}=Fx_t + \varepsilon_Q,
\end{equation}
where 
\begin{align}
	F&=\matr{1&0&1 &0 \\ 0 & 1 & 0 &1 \\ 0 &0 & 1 & 0\\ 0 & 0 & 0 & 1}\\ 
	p(\varepsilon_Q)&=\N\left(\varepsilon_Q\; \middle| \;0,\matr{0&0&0 &0 \\ 0 & 0 & 0 &0 \\ 0 &0 & 0.04 & 0\\ 0 & 0 & 0 & 0.04}\right).
\end{align}
The prior is normal with covariance 
\begin{equation}
P_0=\matr{12&0&0 &0 \\ 0 & 12 & 0 &0 \\ 0 &0 & 1 & 0\\ 0 & 0 & 0 & 1}.
\end{equation}
To have variation in the initial linearizations we sample the initial mean $\mu_0$ from the normal distribution that has a zero mean and covariance  $P_0$. In this test we used Algorithm~\ref{algo:KLPUKF} with \ac{KLD} threshold $\nonllimit =0$, so that each measurement is processed separately at every time step. We tested the \ac{KLPUKF} with moments computed using \ac{nEKF2}, \ac{UKF} and \ac{S2KF} \cite{steinbring2013s}. \ac{S2KF} allows to select the number of sigma-points. We used in our tests 64 sigma-points for \ac{S2KF} while the \ac{UKF} used 9 sigma-points and \ac{nEKF2} used 16 sigma-points. The nonlinearity threshold is set to $\nonllimit=0$ so each measurement element is applied separately. We also tested how the estimation accuracy changes if the measurement elements are applied sequentially in a random order. Estimates are also computed with \ac{EKF}.   Because \ac{EKF} produces always a zero nonlinearity the \ac{KLPUKF} cannot improve the result from the applying the measurements in a random order. The routes were simulated 10000 times.

Table~\ref{tbl:results} shows the mean of errors after the first update. The ``all'' column uses all measurements at once, the ``sequential'' column uses the measurements in random order and the third column uses \ac{KLPUKF}. The table shows how the use of \ac{KLPUKF} instead of the standard filter improves the estimate accuracy more than 15\% at the first time step. Also we can see that the application of the measurements sequentially improves the estimation accuracy.
\begin{table}
\caption{Mean of filtering errors at the first time step}
\label{tbl:results}
\centering
\begin{tabular}{r|c|cc|cc}
& All & \multicolumn{2}{c|}{Sequential} &\multicolumn{2}{c}{\ac{KLPUKF}} \\
\ac{EKF} & 2.35 &  1.95 & -16.7\%  \\ 
\ac{nEKF2} & 1.87 &  1.72 & -8.0\% & 1.48 & -20.9\% \\ 
\ac{UKF} &  1.87 &  1.69 & -10.0\% & 1.54 & -17.9\% \\ 
\ac{S2KF}& 1.82 &  1.69 & -7.5\% & 1.46 & -20.0\% \\ 
\end{tabular}
\end{table}

Table~\ref{tbl:results2} shows the results of the last time step of the routes. In these results we can see that the sequential application of measurements improves only the estimation accuracy of \ac{EKF}; other filters have worse accuracy, but the \ac{KLPUKF} improves the accuracy. The smaller improvement gained by the \ac{KLPUKF} can be explained with the observation that the route moved usually outside the sources of the range measurements and the measurement geometry is worse i.e.\ there is less information in the actual measurements, while they are more linear than in the first step.
\begin{table}
\caption{Mean of filtering errors at the last time step}
\label{tbl:results2}
\centering
\begin{tabular}{r|c|cc|cc}
& All & \multicolumn{2}{c|}{Sequential} &\multicolumn{2}{c}{\ac{KLPUKF}} \\ \hline
\ac{EKF} &  2.17 &  1.99 & -8.1\% & \\ 
\ac{nEKF2} & 1.59 &  1.64 & 3.1\% & 1.50 & -5.6\% \\ 
\ac{UKF} &  1.67 &  1.70 & 1.5\% & 1.57 & -6.2\% \\ 
\ac{S2KF}&  1.59 &  1.62 & 1.5\% & 1.49 & -6.5\% \\ 
\end{tabular}
\end{table}

This example showed that the use of \ac{KLPUKF} gives accuracy improvements and that the use of independent measurements sequentially does not necessarily improve the accuracy.

To get an insight of how to choose $\nonllimit$ we evaluated the previous example with the \ac{S2KF}, which we assume to be the best approximation of the \ac{GGF} of used filters, with different values for $\nonllimit$. Values were set to range from 0 to 2 with interval 0.1 and with $\nonllimit=\infty$, which corresponds to the standard \ac{S2KF}.

Figure~\ref{fig:kldiv} shows how the mean error varies with $\nonllimit$. First of all we can see that in this test setup all measurements had $\eta < 1.5$ so having  $\nonllimit\geq 1.5$ does not change the outcome at all. At the last time step the saturation point of error is already achieved at $\nonllimit\geq 1.5$. At the first and last time step having  $\nonllimit=0.2$ would not increase the error much. At the beginning of the track  $\nonllimit=0.5$ has only a small effect compared to having $\nonllimit=0.5$ at the end of the track. This behavior may be caused by the accumulation of the errors during the track.

\begin{figure}
\centering
	\includegraphics[width=0.8\columnwidth,clip=true,trim=0cm 0cm 0cm 0cm]{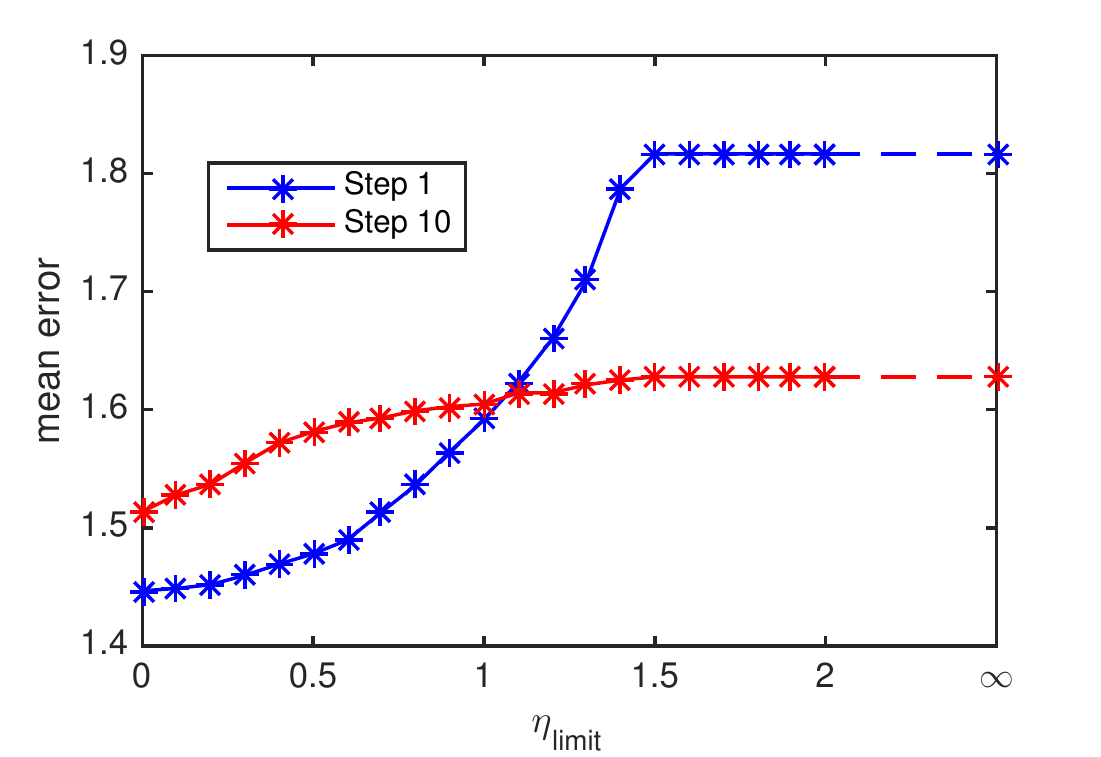}
	\caption{Mean errors as a function of $\nonllimit$}
	\label{fig:kldiv}
\end{figure}

\end{example}

\section{Conclusions}
\label{sec:conclusions}
We have presented \ac{KLPUKF}, an algorithm to perform Bayesian updates with nonlinear multi-dimensional measurements. The key idea is to perform the update in several steps such that the measurement elements with lower nonlinearity are processed first. The measurement elements' processing order is selected so that the joint \ac{KLD} of the state and the measurement is minimized. 

The proposed algorithm can be used with \ac{GGF} approximations, such as sigma-point filters, to improve the estimation accuracy as the nonlinearity measure (\ac{KLD}) can be computed using variables that are already computed within the filter. Simulation results have demonstrated better accuracy of \ac{KLPUKF} compared to several existing filters.

\section*{Acknowledgement}
M. Raitoharju works in OpenKin project that is funded by the Academy of Finland.

\bibliographystyle{IEEEtran}
\bibliography{references}
\balance

\appendix

\begin{proof} \ac{GGF} estimate does not change when a linear transformation is applied to measurements
\label{proof:estc}

Substituting \eqref{equ:transmeas} to the integrals \eqref{equ:i1}-\eqref{equ:i3} the new expectations are 
\begin{align}
	\tilde{\predmeas} &= \int \transmat \measfun(\state)p(\state) \d \state = \transmat \predmeas \label{equ:ii1}\\
	\tilde{\jacobian} &= \int \left( \state - \priormean\right) \left( \transmat\measfun(\state)  - \transmat\predmeas\right)^Tp(\state) \d \state = \jacobian \transmat^T  \label{equ:ii2}\\
	\tilde{\meascov} &= \int \left( \transmat\measfun(\state)  - \transmat\predmeas\right)\left(\transmat \measfun(\state)  -\transmat \predmeas\right)^Tp(\state)  \d \state = \transmat\meascov\transmat \label{equ:ii3}
\end{align}
and the measurement noise covariance of the transformed measurement is
\begin{align}
\tilde{\noisecov}=\transmat \noisecov \transmat^T. \label{equ:noisecov2}
\end{align}
Substituting these into \eqref{equ:mean}, \eqref{equ:posteriorcov}, \eqref{equ:innocov}, and \eqref{equ:kgain}
\begin{align}
	\tilde \innocov & =D(\meascov + \noisecov)D^T \\
	\tilde \kalmangain &= \jacobian D^T D^{-T} \innocov^{-1} D^{-1} =\kalmangain D^{-1} \\
	\tilde \posteriormean & = \priormean + \kalmangain D^{-1}(D \meas - D \predmeas)= \posteriormean\\
	\tilde \posteriorcov & = \priorcov - \kalmangain D^{-1}D\innocov D^T D^{-T}D^T\kalmangain^T= \posteriorcov
\end{align}
we see that the posterior computed using the transformed measurement does not change.
\end{proof}

\begin{proof} Total \ac{KLD} nonlinearity does not change under linear transformation 

\label{proof:nonldoesnotchange}
Using \eqref{equ:ii2} and \eqref{equ:ii3} in \eqref{equ:nonn} we get  
\begin{align}
\begin{aligned}
	\tilde \nonlvar & = \tilde \meascov - \tilde \jacobian^T \left(\priorcov\right)^{-1} \tilde \jacobian  \\
	&= D( \meascov - \jacobian^T \left(\priorcov\right)^{-1} \tilde \jacobian) D^T \\
	&= D \nonlvar D^T
\end{aligned}
\end{align}
By substituting this into \eqref{equ:nonl} with transformed measurement noise covariance  \eqref{equ:noisecov2}
\begin{align}
\begin{aligned}
\tilde{\nonl} &=\frac{1}{2}\log \left| I + (\transmat\noisecov \transmat^T)^{-1}\transmat\nonlvar \transmat^T  \right| \\
&= \frac{1}{2}\log\left| \transmat^{-T}\transmat^T + \transmat^{-T} \noisecov^{-1}\nonlvar\transmat^{-T} \right| \\
&= \frac{1}{2}\log\left(\left| \transmat^{-T}\right|\left| I  + \noisecov^{-1}\nonlvar\right|\left| \transmat^{T}\right|\right) \\
&= \frac{1}{2}\log \left| I + \noisecov^{-1}\nonlvar  \right| \\
&= \nonl
\end{aligned}
\end{align}
we see that the total nonlinearity does not change.
\end{proof}

\begin{proof} Nonlinearity of an element is minimized

\label{proof:maxnl}
The proof is very similar to one in the Appendix~B in \cite{PUKF}.
When  $\tilde R=I$ and $\tilde \nonlvar$ are diagonal the total nonlinearity is 
\begin{align}
	\tilde{\nonl} &= \frac{1}{2}\log | I+ \tilde{\nonlvar}| \\ &= \sum_{i=1}^\measdim \frac{1}{2}\log( 1 + \tilde{\nonlvar}_{[i,i]}) \label{equ:nonlvaluetot}
\end{align}
and the nonlinearity corresponding to $i$th measurement element is
\begin{align}
	\tilde{\nonl}_i &=\frac{1}{2}\log( 1 + \tilde{\nonlvar}_{[i,i]}) \label{equ:nonlvalue}
\end{align}

We will show that the smallest diagonal element of $\tilde{\Upsilon}$ is as small as possible under a linear transformation that preserves  $\noisecov=I$ and further that the second smallest diagonal element is as small as possible, when the next smallest is as small as possible etc. If the measurement model is transformed by multiplying it with matrix $V$, the transformed variables are $\hat \nonlvar= V \tilde \nonlvar V^T$ and $\hat \noisecov=VIV^T=VV^T$. Because we want to have $\noisecov=I$, $V$ has to be unitary. The $i$th diagonal element of the transformed matrix $\tilde \nonlvar$ is $v_i^T \tilde{\nonlvar}  v_i = \sum_{j=1}^\measdim v_{i,[j]}^2 \tilde{\nonlvar} \el{j,j}$, where $v_i$ is the $i$th column of $V$. Because $V$ is unitary $\sum_{j=1}^\measdim v_{i,[j]}^2=1$ and the $i$th diagonal element of the transformed matrix $\hat\nonlvar$ is
\begin{equation}
\sum_{j=1}^d v_{i, [j]}^2 \tilde{\nonlvar}_{[j,j]} \geq \sum_{j=1}^d v_{i, [j]}^2 \min_j\{\tilde{\nonlvar}_{[j,j]}\} = \min_j \{{\nonlvar}_{[j,j]}\}.
\end{equation}
Thus, the new diagonal element cannot be smaller than the smallest diagonal element of ${\tilde{\nonlvar}}$.

If the smallest element is the first element of the diagonal the possible transformation for the second smallest element is
\begin{equation}
	\hat{\nonlvar} = \matr{1 & 0^T \\ 0  & V} \tilde{\nonlvar} \matr{1 & 0^T \\ 0  & V^T} .
\end{equation}
With the same reasoning as given already the second diagonal has to be already the smallest possible. Inductively this applies to all diagonal elements. One could also show that now the measurement element corresponding to the maximal element of $\tilde{\nonlvar}$ has the largest possible \ac{KLD}.
\end{proof}

\end{document}